\documentclass[11pt]{article}
\usepackage{enumerate}
\usepackage{amssymb,a4wide,latexsym,makeidx,epsfig,fleqn}
\usepackage{amsthm}
\usepackage{amsmath}
\usepackage{enumerate}
\newtheorem{thm}{Theorem}[section]

\newtheorem{cor}[thm]{Corollary}
\newtheorem{lem}[thm]{Lemma}
\newtheorem{defn}[thm]{Definition}

\begin{document}
\textwidth 150mm \textheight 225mm
\title{No signed graph with the nullity $\eta(G,\sigma)=|V(G)|-2m(G)+2c(G)-1$ \footnote{This work is supported by the National Natural Science Foundations of China (No. 11901253), the Natural Science Foundation for Colleges and Universities in Jiangsu Province of China (No. 19KJB110009), and the Science Foundation of Jiangsu Normal University (
No. 18XLRX021).}}
\author{{Yong Lu\footnote{Corresponding author.}, Jingwen Wu}\\
{\small  School of Mathematics and Statistics, Jiangsu Normal University,}\\ {\small  Xuzhou, Jiangsu 221116,
People's Republic
of China.}\\
{\small E-mail: luyong@jsnu.edu.com}}

\date{}
\maketitle
\begin{center}
\begin{minipage}{120mm}
\vskip 0.3cm
\begin{center}
{\small {\bf Abstract}}
\end{center}
{\small Let $G^{\sigma}=(G,\sigma)$ be a signed graph and $A(G,\sigma)$ be its adjacency matrix. Denote by $m(G)$ the matching number of $G$. Let $\eta(G,\sigma)$ be the nullity of $(G,\sigma)$. He et al. [Bounds for the matching number and cyclomatic number of a signed graph in terms of rank, Linear Algebra Appl. 572 (2019), 273--291] proved that
$$|V(G)|-2m(G)-c(G)\leq\eta(G,\sigma)\leq |V(G)|-2m(G)+2c(G),$$ where $c(G)$ is the dimension of cycle space of $G$. Signed graphs reaching the lower bound or the upper bound are respectively characterized by the same paper.
In this paper, we will prove that no signed graphs with nullity $|V(G)|-2m(G)+2c(G)-1$. We also prove that there are  infinite signed  graphs with nullity $|V(G)|-2m(G)+2c(G)-s,~(0\leq s\leq3c(G), s\neq1)$ for a given $c(G)$.

 }
\end{minipage}
\end{center}

\section{Introduction }



Let $G=(V(G),E(G))$ be a simple graph, where $V(G)$ is the vertex
set of $G$ and $E(G)$ is the edge set of $G$. Let $V(G)=\{v_{1},v_{2},\ldots,v_{n}\}$. Then the \emph{adjacency matrix} $A(G)$ of $G$ is the symmetric $n\times n$ matrix with entries $A(i,j)=1$ (or written as $a_{ij}=1$) if and only if $v_{i}v_{j}\in E(G)$ and zeros elsewhere. The \emph{nullity}(resp. \emph{rank}) of $G$, denoted by $\eta(G)$(resp. $r(G)$), is the multiplicity of 0 (resp. non-zero) eigenvalue of $A(G)$. A \emph{pendant vertex} is defined as a vertex with degree one. The unique neighbour of the pendant vertex is called \emph{quasi-pendant vertex}. If a cycle of $G$ contains only a vertex of degree three, and the remaining vertices are of degree two, then the cycle is called a \emph{pendant cycle}.

Assume that $V_{1}\subseteq V(G)$ and $V_{1}\neq \phi$, Let $G[V_{1}]$ be the \emph{induced subgraph} of $G$, whose vertex set is $V_{1}$ and edge set is the set of those edges of $G$ that have both ends in $V_{1}$. Denoted by $G-V_{1}$ the induced subgraph obtained from $G$ by deleting each vertex in $V_{1}$ and its incident edges. For convenience, if $V_{1}=\{x\}$, we write $G-x$ for $G-\{x\}$. For the induced subgraph $H_{1}$ of $G$, $H_{1}+x$ is defined as the subgraph of $G$ induced by the vertex set $V(H_{1})\cup \{x\}$.

The \emph{matching number} of $G$, denoted by $m(G)$, is the size of a maximum matching of $G$.
Let $M$ be a maximum matching of $G$, if there exists the edge $e\in M$ such that the vertex $v\in e$, $v$ is called \emph{$M$-saturated}, otherwise, $v$ is called \emph{$M$-unsaturated}.  An \emph{$M$-alternating path} of $G$ is defined as a path whose edges are alternately in the edge sets $E\backslash M$ and $M$. An \emph{$M$-augmenting path} is defined as a path whose starting vertex and ending vertex are $M$-unsaturated.

The \emph{distance} between $u$ and $v$, denote by $d(u,v)$, is the length of the shortest path from the vertex $u$ to $v$. The length of the shortest cycle of $G$ is called the \emph{girth} of $G$, denote by $g(G)$.  Let $c(G)$ be \emph{the dimension of cycle space} of $G$, $c(G)=|E(G)|-|V(G)|+\theta(G)$, where $\theta(G)$ is the number of connected components of $G$. If the  cycles (if any) of $G$ are pairwise vertex-disjoint, the acyclic graph $T_{G}$ is obtained from  contracting each cycle of $G$ into a vertex, which is called \emph{cyclic vertex}. Let $W_{G}$ (resp. $U$) be the vertex set consisting of all cyclic vertices (resp. all non-cyclic vertices) in $T_{G}$, $V(T_{G})=W_{G}\cup U$. Furthermore, the graph is obtained from $T_{G}$ by deleting all cyclic vertices, denoted by $[T_{G}]$.

A \emph{signed graph} $G^{\sigma}=(G, \sigma)$ consists of a simple graph $G$ with edge set $E$ and a mapping $\sigma: E\rightarrow\{+, -\}$. $G$ is called the \emph{underlying graph} of $(G, \sigma)$. The \emph{adjacency matrix}  of $(G, \sigma)$, denoted by $A(G, \sigma)=a^{\sigma}_{ij}=\sigma(v_{i}v_{j})a_{ij}$, where $a_{ij}\in A(G)$. We use $r(G, \sigma)$ to denote the \emph{rank} of a signed graph $(G, \sigma)$.

Denote by  $C^{\sigma}_{n}$ a signed cycle of order $n$. The  \emph{sign} $sgn(C^{\sigma}_{n})$ of $C^{\sigma}_{n}$ is defined as $\prod_{e\in E(C^{\sigma}_{n})}\sigma(e)$. If $sgn(C^{\sigma}_{n})=+$ (or $sgn(C^{\sigma}_{n})=-$), then $C^{\sigma}_{n}$ is said to be \emph{positive} (or \emph{negative}). If all the cycles  of $(G,\sigma)$ are positive, then $(G,\sigma)$ is \emph{balanced}, and \emph{unbalanced} otherwise.

\begin{figure}[htbp]
 \centering
  \includegraphics[scale=0.65]{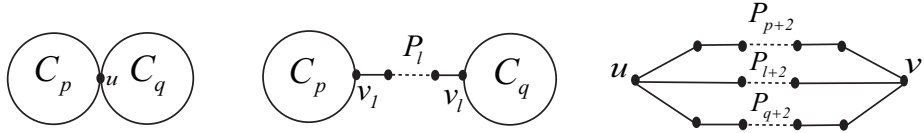}
 \caption{$\infty(p,1,q)$, $\infty(p,l,q)$ and $\theta(p,l,q)$.}
\end{figure}

A connected graph $G$ is called \emph{bicyclic} if $|E(G)|=|V(G)|+1$. 
The \emph{base} of a bicyclic graph $G$ is the unique connected bicyclic subgraph of $G$ containing no pendant vertices.

Let $\infty(p,l,q)$ (see Figure 1) be the graph obtained from $C_{p}$, $C_{q}$ and $P_{l}$ by identifying $v$ with $v_{1}$ and $u$ with $v_{l}$, respectively.
 A bicyclic graph $G$ is called an $\infty$-graph, if $G$ contains
$\infty(p,l,q)$ as its base.
Let $P_{p+2}$,~$P_{q+2}$,~$P_{l+2}$ be three paths, where $min\{p,l,q\}\geq0$ and at most one of $p,~l,~q$ is 0. Let $\theta(p,l,q)$ (see Figure 1) be the graph obtained from $P_{p+2}$,~$P_{l+2}$ and $P_{q+2}$ by identifying the three initial vertices and terminal vertices. The bicyclic graph containing
$\theta(p,l,q)$ as its base is called a $\theta$-graph.

A connected graph $G$ is called \emph{tricyclic} if $|E(G)|=|V(G)|+2$.
As shown in Figure 2, there are eight types of bases of tricyclic graph, denoted by $T_{i}$, $i=1,\ldots,8$. The tricyclic graph $G$ can be obtained from $T_{i}$ by attaching trees to its vertices. 

\begin{figure}[htbp]
 \centering
  \includegraphics[scale=0.5]{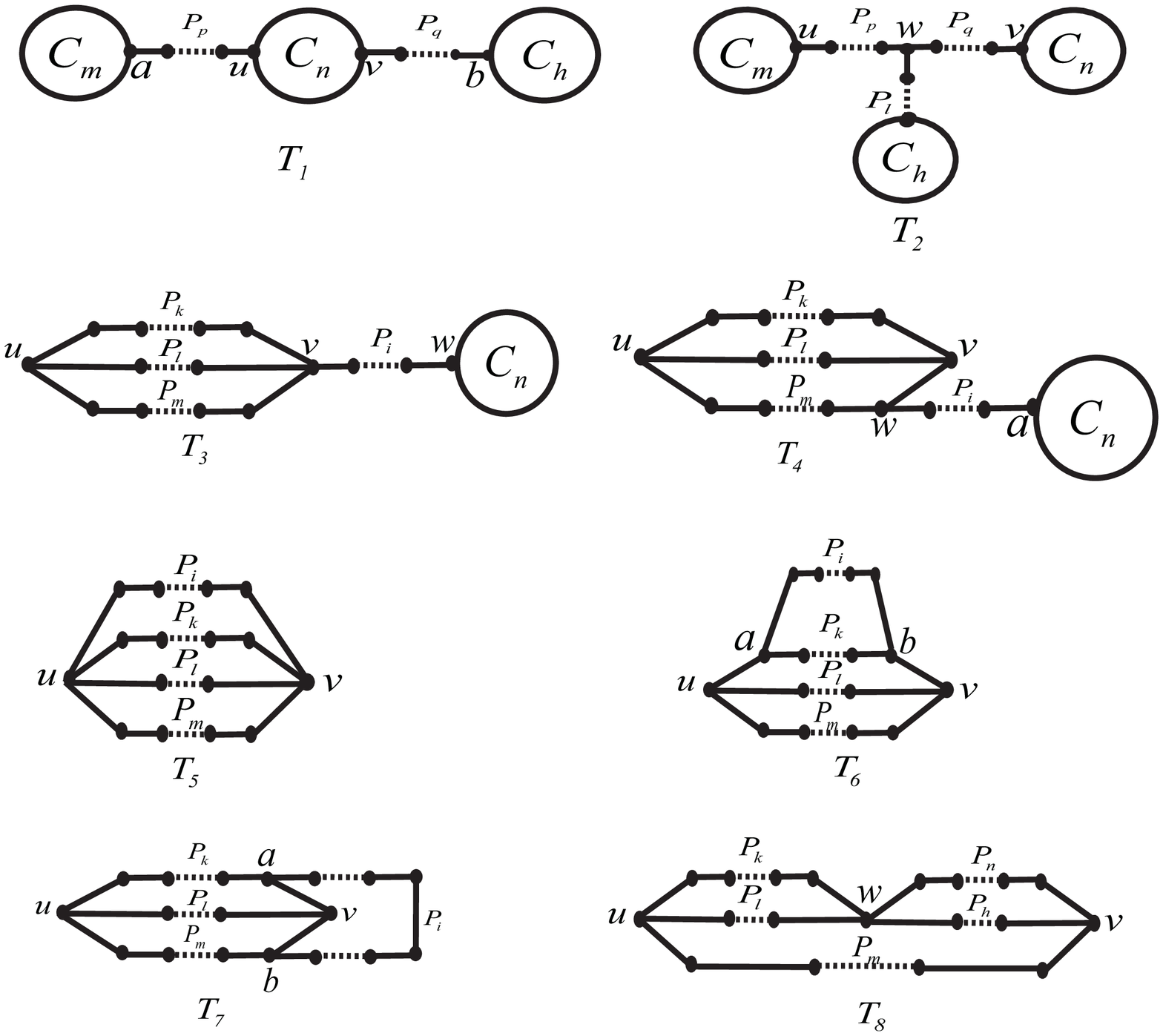}
 \caption{$T_{1}$--$T_{8}$}
\end{figure}

Collatz et al. \cite{CS}  had wanted to obtain all graphs of order $n$ with $\eta(G)>0$. Until today, this problem is also unsolved.
In recent years, the investigation on the nullity (or rank) of simple graphs is an significant topic in the theory of graph spectra. In 2014, Wang and Wong
\cite{WLWD} investigated the relationship between the nullity of the simple graph and the matching number, and proved that $|V(G)|-2m(G)-c(G)\leq \eta(G)\leq |V(G)|-2m(G)+2c(G)$. In 2015, Song, Song and Tam \cite {SST} characterized the graphs with the nullity $\eta(G)= |V(G)|-2m(G)+2c(G)$. In 2016, Wang \cite {WL} characterized the graphs with nullity $\eta(G)=|V(G)|-2m(G)-c(G)$. In 2019, Li and Guo \cite{LXGJM} proved that there is no graph with nullity $\eta(G)=|V(G)|-2m(G)+2c(G)-1$, and for fixed $c(G)$, infinitely many connected graphs with nullity $\eta(G)=|V(G)|-2m(G)+2c(G)-s$, where $0\leq s\leq 3c(G)$ and $s\neq1$, are also constructed.

In 2018, Lu, Wang and Zhou \cite{LYWLG} proved that $\eta(G)-2c(G)\leq \eta(G,\sigma) \leq \eta(G)+2c(G)$ for an unbalanced signed graph with $c(G)\geq1$, and all corresponding extremal graphs are characterized.
In 2019, He, Hao and Lai \cite {HSJHRX} proved that for every signed graph, $|V(G)|-2m(G)-c(G)\leq \eta(G,\sigma)\leq |V(G)|-2m(G)+2c(G)$. Moreover, signed graphs reaching the lower bound or the upper bound are respectively characterized in \cite{HSJHRX}. In this paper, we prove that no signed graph with nullity
$\eta(G,\sigma)=|V(G)|-2m(G)+2c(G)-1$, and for fixed $c(G)$, infinitely many connected signed graphs with nullity $\eta(G,\sigma)=|V(G)|-2m(G)+2c(G)-s$, where $0\leq s\leq 3c(G)$ and $s\neq1$, are also constructed.

Our main results are the following three theorems.

\begin{thm}\label{thm:1.1}

Suppose $(G,\sigma)$ is an unbalanced signed unicyclic graph with an unique cylce $C_{q}^{\sigma}$, then
\begin{enumerate}[(1)]
\item  $\eta(G,\sigma)=|V(G)|-2m(G)-1$, if $q\equiv1(mod~2)$,~$m([T_{G}])=m(T_{G})$~(Lemma 4.3~\cite{HSJHRX});
\item  $\eta(G,\sigma)=|V(G)|-2m(G)+2$, if $q\equiv2(mod~4)$,~$m([T_{G}])=m(T_{G})$~(Lemma 5.3~\cite{HSJHRX});
\item  $\eta(G,\sigma)=|V(G)|-2m(G)$, otherwise.
\end{enumerate}
\end{thm}

\begin{thm}\label{thm:1.2}

For any unbalanced signed graph $(G,\sigma)$, $$\eta(G,\sigma)\neq|V(G)|-2m(G)+2c(G)-1.$$
\end{thm}

\begin{thm}\label{thm:1.3}
Let $(G,\sigma)$ be an unbalanced signed graph. If $c(G)$ is fixed, then there exists infinitely connected unbalanced signed graphs, such that $\eta(G,\sigma)=|V(G)|-2m(G)+2c(G)-s$, where $0 \leq s \leq 3c(G), s\neq1$.

\end{thm}

\section{Preliminaries}

\begin{lem} \label{le:2.13}\cite{Bondy}
A matching $M$ of $G$ is a maximum matching if and only if $G$ contains no M-augmenting path.
\end{lem}
\begin{lem} \label{le:2.8}\cite{MXBWD}
Let $G$ be a graph with a pendant vertex $u$, and $v$ is the unique quasi-pendant vertex of $u$. Then $m(G)=1+m(G-v)=1+m(G-u-v)$.

\end{lem}
\begin{lem} \label{le:2.9}\cite{CCJH}
Let $G$ be a simple graph. Then $m(G)-1\leq m(G-v) \leq m(G)$, for any vertex $v \in V(G)$.
\end{lem}

\begin{lem} \label{le:2.10}\cite{MXBWD}
Let $G$ be a graph obtained by joining a vertex of an even cycle $C$ by an edge to a vertex y of a connected graph $H$. Then $m(G)=m(C)+m(H)$.
\end{lem}

\begin{lem} \label{le:2.14}\cite{LXGJM}
For any undirected graph $G$, $\eta(G)\neq|V(G)|-2m(G)+2c(G)-1.$
\end{lem}

\begin{lem} \label{le:2.12}\cite{WDMXB}
Let $G$ be a graph with $x \in V(G)$,

\begin{enumerate}[(1)]
  \item $c(G-x)=c(G)$  if $x$ lies outside any cycle of $G$;
  \item $c(G-x)\leq c(G)-1$  if $x$ lies on a cycle of $G$;
  \item $c(G-x)\leq c(G)-2$  if $x$ is a common vertex of distinct cycles of $G$;
  \item if the cycles of $G$ are pairwise vertex-disjoint, then $c(G)$ precisely equals the number of cycles in $G$.
\end{enumerate}

\end{lem}

\begin{lem} \label{le:2.1}\cite{BJBKDG}
Let $x$ be a vertex of $(G,\sigma)$. Then $\eta(G,\sigma)-1 \leq \eta(G-x,\sigma) \leq \eta(G,\sigma)+1$.
\end{lem}

\begin{lem} \label{le:2.4}\cite{FYZWY}
Let $(G,\sigma)$ be a signed graph containing a pendant vertex $u$, $uv\in E(G)$ and $(H,\sigma)=(G,\sigma)-u-v$. Then $\eta(G,\sigma)=\eta(H,\sigma)$.
\end{lem}

\begin{lem} \label{le:2.2}\cite{HSJHRX}
Let $(G,\sigma)$ be a connected signed graph. Then $|V(G)|-2m(G)-c(G) \leq \eta(G,\sigma) \leq |V(G)|-2m(G)+2c(G)$.
\end{lem}

\begin{lem} \label{le:2.6}\cite{FYZWY}
Let $P_{n}^{\sigma}$ be a signed path. Then,
\begin{equation}
\eta(P_{n}^{\sigma})=\left\{\begin{array}{ll}
0, ~if~n\equiv0(mod~2);\\\\
1, ~if~n\equiv1(mod~2).
\end{array}
\right.\nonumber
\end{equation}
\end{lem}

\begin{lem} \label{le:2.7}\cite{FYZWY}
Let $C^{\sigma}_{p}$ be a signed cycle. Then,

\begin{equation}
\eta(C^{\sigma}_{p})=\left\{\begin{array}{ll}
0, ~if~p\equiv1(mod~2);\\\\
0, ~if~p\equiv0(mod~4)~and~\sigma(C_{p})=-;\\\\
0, ~if~p\equiv2(mod~4)~and~\sigma(C_{p})=+;\\\\
2, ~if~p\equiv0(mod~4)~and~\sigma(C_{p})=+;\\\\
2, ~if~p\equiv2(mod~4)~and~\sigma(C_{p})=-.
\end{array}
\right.\nonumber
\end{equation}
\end{lem}

\begin{lem} \label{le:2.5}\cite{HYPLJS}
Let $(G,\sigma)$ be a signed graph. Then $(G,\sigma)$ is balanced if and only if $(G,\sigma)\thicksim (G,+)$.

\end{lem}

\begin{lem} \label{le:2.11}\cite{LYWLG}
Let $(G,\sigma)$ be a signed graph, if $(G_{1},\sigma)$,$(G_{2},\sigma)$,$\cdots$,$(G_{t},\sigma)$ are all the connected components of $(G,\sigma)$. Then $\eta(G,\sigma)=\sum_{i=1}^{t} \eta(G_{i},\sigma)$.

\end{lem}

\begin{lem} \label{le:2.3}\cite{HSJHRX}
Let $(G,\sigma)$ be a connnected signed graph. Then $\eta(G,\sigma)=|V(G)|-2m(G)+2c(G)$ if and only if all the following conditions hold for $(G,\sigma)$:
\begin{enumerate}[(1)]
\item  the cycles (if any) of $(G,\sigma)$ are pairwise vertex-disjoint.
\item  for each cycle (if any) $C_{q}^{\sigma}$ of $(G,\sigma)$, either $q\equiv0(mod~4)$ and $\sigma(C_{q})=+$ or $q\equiv 2(mod~4)$ and $\sigma(C_{q})=-$.
\item  $m(T_{G})=m([T_{G}])$.
\end{enumerate}
\end{lem}


\section{Proof of Theorem 1.1}
In this section, we  will proved that there is no unbalanced signed unicyclic graph with the nullity $\eta (G,\sigma)=|V(G)|-2m(G)+2c(G)-1$. In order to proof Theorem \ref{thm:1.1}, we need to introduce some definitions and lemmas.

In simple graph, we introduced the concept of elementary graph, either the complete graph with order two $K_{2}$ or the cycle $C_{q},~q\geq3$. If the graph whose each component is elementary graph, the graph is called the basic graph. The set of signed graphs, denoted by $\mathbf{U_{i}}$, is the set of all basic graphs with $i$ vertices in $(G,\sigma)$. For any basic graph $(u,\sigma)\in \mathbf{U_{i}}$, $p(u), ~c(u)$ and $s(u)$ are defined as the number of components, the number of cycles and the number of the negative edges in the cycle of $(u,\sigma)$, respectively.

\begin{defn}\label{de:3.2}
Let $G$ be an unicyclic graph with an unique cycle $C_{q}$. Then,

\begin{enumerate}[(1)]
  \item $E_{1}$: The edge set of $G$ between  $C_{q}$ and $[T_{G}]$.
  \item $F_{1}$: The matching set of $G$ with $m(G)$ edges.
  \item $F_{2}$: The matching set of $[T_{G}]$ with $m([T_{G}])$ edges.
  \item $F_{1}^{'}$: The matching set of $G$ with $m(G)$ edges, each of which has at least one edge in $E_{1}$.
  \item $F_{1}^{''}$: The matching set of $G$ with $m(G)$ edges, and $M\cap E_{1}=\phi$ for all $M\in F_{1}$.
\end{enumerate}
\end{defn}

By Definition \ref{de:3.2}, we have $F_{1}=F_{1}^{'}\cup F_{1}^{''}$. We can obtain the following corollary:

\begin{cor}\label{co:1}
 Let $C_{q}$ be an even cycle. Then,
\begin{enumerate}[(1)]
  \item If $F_{1}^{'}=\phi$, the maximum matching of $G$ is the union of a maximum matching of $C_{q}$ and a maximum matching of $G-C_{q}$. Then  $|F_{1}|=|F_{1}^{''}|=2|F_{2}|$.
  \item If $F_{1}^{'}\neq \phi$, then $|F_{1}|=|F_{1}^{'}|+|F_{1}^{''}|> 2|F_{2}|$.

\end{enumerate}
\end{cor}

\begin{lem} \label{le:3.1}\cite{LYYLH}
Let $(G,\sigma)$ be a signed graph on $n$ vertices and $P_{(G,\sigma)}(\lambda)=\Sigma_{i=0}^{n}a_{i} \cdot \lambda^{n-i}$ be a characteristic polynomial of $A(G,\sigma)$. Then $a_{i}=\sum_{(u,\sigma) \in \mathbf{U_{i}}}(-1)^{p(u)+s(u)}\cdot 2^{c(u)}$, $i \in \{1,2,\ldots,n\}$.

\end{lem}

\begin{lem} \label{le:3.2}\cite{HSJHRX}
Let $(G,\sigma)$ be a connected unbalanced signed unicyclic graph whose unique cycle is $C_{q}^{\sigma}$. Then $\eta(G,\sigma)=|V(G)|-2m(G)+2$ if and only if $q\equiv 2(mod~4)$ and $m(T_{G})=m([T_{G}])$.

\end{lem}
\begin{lem} \label{le:3.3}\cite{HSJHRX}
Let $(G,\sigma)$ be an unicyclic signed graph which contains the unique signed cycle $C_{q}^{\sigma}$. Then $\eta(G,\sigma)=|V(G)|-2m(G)-1$ if and only if $q$ is odd and $m(T_{G})=m([T_{G}])$.

\end{lem}

\begin{lem} \label{le:3.4}\cite{RSZY}
Let $G$ be a graph with at least one cycle. Suppose that all cycles of $G$ are pairwise vertex-disjoint and each cycle is odd, then $m(T_{G})=m([T_{G}])$ if and only if $m(G)=\sum_{C\in \mathcal{L}(G)}m(C)+m([T_{G}])$, where $\mathcal{L}(G)$ denotes the set of all cycles in $G$.

\end{lem}

By Lemma \ref{le:3.4}, we can obtain the following corollary.
\begin{cor}\label{co:2}
Let $G$ be an unicyclic graph with the unique cycle $C_{q}$, $q$ is odd. Then
 $m(T_{G})= m([T_{G}])$ if and only if $m(G)=m([T_{G}])+m(C_{q})$.

\end{cor}

Next, for the unicyclic graph with an even cycle, similar to Lemma \ref{le:3.4}, we give the following lemma.

\begin{lem} \label{le:3.5}
Let $G$ be an unicyclic graph with a cycle $C_{q}$, $q$ is even. Then $m(T_{G})=m([T_{G}])$ if and only if $m(G)=m(C_{q})+m([T_{G}])$ and $M\cap E_{1}=\phi$ for all $M\in F_{1}$.

\end{lem}

\noindent\textbf{Proof.} Suppose that $G-C_{q}$ has $p$ connected components, says $S_{1},S_{2},\ldots,S_{p}$. Since $[T_{G}]=G-C_{q}$, then $S_{1},S_{2},\ldots,S_{p}$ are also all connected components of $[T_{G}]$. Assume that the vertex $v_{c}$ of $T_{G}$ is correspond to the cycle $C_{q}$ of $G$, which implies that $W_{G}=\{v_{c}\}$ in $T_{G}$.

In $T_{G}$, let $\widetilde{S_{k}}$ ($k=1,\ldots,p$) be the subgraph induced by $V(S_{k})$ and all vertices adjacent to it. Similarly, in $G$, let $\widehat{S_{k}}$ be the subgraph induced by $V(S_{k})$ and all vertices adjacent to it. Since $G$ is an unicyclic graph, then there exists at most an edge between $V(S_{k})$ and $V(C_{q})$, 
thus for all $k=1,\ldots,p$, the subgraph $\widetilde{S_{k}}$ and $\widehat{S_{k}}$ are isomorphic.

\textbf{Necessity:}
Note that $m(T_{G})=m([T_{G}])$, and  $m([T_{G}])=\sum_{k=1}^{p}m(S_{k})$. Then each maximum matching of $[T_{G}]$ is the maximum matching of $T_{G}$, denoted by $M_{T}$ the maximum matching of $[T_{G}]$, thus $|M_{T}|=\sum_{k=1}^{p}m(S_{k})$. We can prove $m(S_{k})=m(\widetilde{S_{k}})$ for $k=1,\ldots,p$. Otherwise, there exists $m(S_{k_{0}})<m(\widetilde{S_{k_{0}}})$, where $k_{0}\in \{1,\ldots,p\}$. For $M_{T}$, replace the maximum matching of $S_{k_{0}}$ to the maximum matching of $\widetilde{S_{k_{0}}}$, we obtained a matching of $T_{G}$, denoted by $M_{T}^{'}$. Then，
$|M_{T}^{'}|= m(\widetilde{S_{k_{0}}})+\sum_{k\neq k_{0}}m(\widetilde{S_{k}}) > m(S_{k_{0}})+\sum_{k\neq k_{0}}m(S_{k})= |M_{T}|$.

This contradicts $|M_{T}|=m(T_{G})\geq |M_{T}^{'}|$. Hence $m(S_{k})=m(\widetilde{S_{k}})$, for all $k=1,\ldots,p$. Since the subgraph $\widetilde{S_{k}}$ and $\widehat{S_{k}}$ are isomorphic, then $m(\widehat{S_{k}})=m(S_{k})=m(\widetilde{S_{k}})$, for all $k=1,\ldots,p$. Let $M_{1}$ be a maximum matching of $G$. Then
$m(G)=|M_{1}|=|M_{1}\cap E(C_{q})|+\sum_{k=1}^{p}|M_{1}\cap E(\widehat{S_{k}})| \leq m(C_{q})+\sum_{k=1}^{p}m(\widehat{S_{k}})=m(C_{q})+\sum_{k=1}^{p}m(S_{k})$. By $m(T_{G})=m([T_{G}])$, we have $m(G)\geq m(C_{q}) + \sum_{k=1}^{p}m(S_{k})$.

Thus $$m(G)= m(C_{q}) + \sum_{k=1}^{p}m(S_{k}),$$

that is, $$m(G)= m(C_{q}) + m([T_{G}]).$$
Next, we need to prove $M\cap E_{1}=\phi$ for any maximum matching $M$ of $G$, $M\in F_{1}$. Suppose on the contrary, there exists $M_{2} \in F_{1}$ such that $ M_{2} \cap E_{1} \neq \phi$. Let $uv \in E_{1} \cap M_{2}$, where $v\in V(C_{q}), u\in V(S_{k_{1}})$, and $k_{1} \in \{1,\ldots,p\}$. Since $C_{q}$ has a perfect matching, which shows that $m(C_{q}-v)=m(C_{q})-1<m(C_{q})$. In $\widehat{S_{k_{1}}}$, $u$ is a quasi-pendant vertex, by Lemma \ref{le:2.8}, $m(S_{k_{1}}-u)=m(\widehat{S_{k_{1}}}-u-v)=m(\widehat{S_{k_{1}}})-1=m(S_{k_{1}})-1$.
Then,
$|M_{2}|=m(C_{q}-v)+(1+m(\widehat{S_{k_{1}}}-u-v))+\sum_{k\neq k_{1}}|M_{2} \cap \widehat{S_{k}}|<m(C_{q})+\sum_{k=1}^{p}m(\widehat{S_{k}})=m(C_{q})+\sum_{k=1}^{p}m(S_{k})=m(G)$. This contradicts $m(G)=|M_{2}|$, thus
$M\cap E_{1}=\phi$ for all $M\in F_{1}$.

\textbf{Sufficiency:}
Note that $m(G)= m(C_{q}) + \sum_{k=1}^{p}m(S_{k})$. We assume that $m(T_{G})\neq m([T_{G}])$, that is, $m(T_{G})=m([T_{G}])+1$ ($G$ is an unicyclic graph). Which implies that each maximum matching of $T_{G}$ must cover the vertex $v_{c}$ of $W_{G}$. There exists a maximum matching $M_{T}^{''}$ of $T_{G}$ such that the edge $wv_{c}\in M_{T}^{''}$, and the vertex $w$ is adjacent to $v_{0}$ the vertex of $C_{q}^{\sigma}$ in $G$.
Since $q$ is even, then $m(C_{q}-v_{0})=m(C_{q})-1$. Assume that $M_{3}$ is a matching of $G$, which is the union of $M_{T}^{''}$ and a maximum matching of $C_{q}-v_{0}$. Then $|M_{3}|=m(T_{G})+m(C_{q}-v_{0})=m([T_{G}])+1+m(C_{q})-1=m(G)$, which implies that $M_{3}$ is a maximum matching of $G$ and $M_{3}\cap E_{1} \neq \phi$, a contradiction. Thus $m(T_{G})= m([T_{G}])$.\quad $\square$

Next, we give Lemmas \ref {le:3.6} and \ref {le:3.7} by applying the characteristic polynomial of the signed graph.
\begin{lem} \label{le:3.6}
Let $(G,\sigma)$ be an unbalanced signed unicyclic graph, which contains the unique cycle $C_{q}^{\sigma}$. If $q$ is odd and $m(T_{G})\neq m([T_{G}])$, then $\eta(G,\sigma)=|V(G)|-2m(G)$.
\end{lem}

\noindent\textbf{Proof.}
Let $|V(G)|=n,~m(G)=m~ and~ m(G-C_{q})=l$. By Lemma \ref{le:3.1}, we have the characteristic polynomial of $(G,\sigma)$, $P_{(G,\sigma)}(\lambda)=\sum_{i=0}^{n}a_{i} \cdot \lambda^{n-i}$,~$a_{i}=\sum_{(u,\sigma) \in \mathbf{U_{i}}}(-1)^{p(u)+s(u)}\cdot 2^{c(u)}$, for any $i \in \{1,2,\ldots,n\}$.
Since $q$ is odd, then $m(C_{q})=\frac{q-1}{2}$, and $C_{q}^{\sigma}$ has not a perfect matching. Note that $m(T_{G})\neq m([T_{G}])$, by Corollary \ref{co:2}, we have $m(G)>m([T_{G}])+m(C_{q})=l+\frac{q-1}{2}$, that is, $2m+1>q+2l$. If $i> 2m$, then $(G,\sigma)$ contains no basic subgraphs with $i$ vertices and $a_{i}=0$. Thus  $P_{(G,\sigma)}(\lambda)=\lambda^{n}+a_{1}\lambda^{n-1}+\cdots +a_{2m}\lambda^{n-2m}=\lambda^{n-2m}(\lambda^{2m}+a_{1}\lambda^{2m-1}+\cdots +a_{2m})$. Therefore
$\eta(G,\sigma)\geq |V(G)|-2m(G)$.

Next, we need to prove $a_{2m}\neq 0$. Note that $a_{2m}=\sum_{(u,\sigma) \in \mathbf{U_{2m}}}(-1)^{p(u)+s(u)}\cdot 2^{c(u)}$, since $C_{q}^{\sigma}$ is an odd cycle, the basic subgraphs of $(G,\sigma)$ with $2m$ vertices contain no cycle. Then $s(u)=c(u)=0,~p(u)=m$, we have $a_{2m}=\sum_{(u,\sigma) \in \mathbf{U_{2m}}}(-1)^{m}\neq 0$. Thus $\eta(G,\sigma)=|V(G)|-2m(G)$.\quad $\square$
\begin{lem} \label{le:3.7}
Let $(G,\sigma)$ be an unbalanced signed unicyclic graph with the unique cycle $C_{q}^{\sigma}$, and $q$ is even. If  $m(T_{G})= m([T_{G}])$ and $q\equiv 0(mod~4)$; or $m(T_{G})\neq m([T_{G}])$, then $\eta(G,\sigma)=|V(G)|-2m(G)$.
\end{lem}

\noindent\textbf{Proof.}
Let $F_{3}=\{u|u=C_{q}\cup M, M \in F_{2}\}$, $m(G)=m$ and $m(G-C_{q})=l$. $C_{q}^{\sigma}$ has odd number of negative edges. Note that $G$ is a bipartite graph. By Lemma \ref{le:3.1}, then the characteristic polynomial of $(G,\sigma)$ can be expressed by $P_{(G,\sigma)}(\lambda)=\sum_{i=0}^{\lfloor \frac{n}{2} \rfloor}b_{i} \cdot \lambda^{n-2i}$, $b_{i}=\sum_{(u,\sigma) \in \mathbf{U_{2i}}}(-1)^{p(u)+s(u)}\cdot 2^{c(u)}$, for any $i \in \{1,2,\ldots,\lfloor \frac{n}{2} \rfloor\}$. 
By the similar proof as in Lemma \ref{le:3.6}, if $i>m$, then $(G,\sigma)$ contains no basic subgraphs with $2i$ vertices and $b_{i}=0$. Hence $$P_{(G,\sigma)}(\lambda)=\lambda^{n}+b_{1}\lambda^{n-2}+\cdots +b_{m}\lambda^{n-2m}=\lambda^{n-2m}(\lambda^{2m}+b_{1}\lambda^{2m-2}+\cdots +b_{m}).$$ Therefore $\eta(G,\sigma)\geq |V(G)|-2m(G)$.

Next, in order to get the result $\eta(G,\sigma)=|V(G)|-2m(G)$, we need to prove $b_{m}\neq 0$.

\textbf{(1).} $m(T_{G})= m([T_{G}])$, and $q\equiv 0(mod~4)$.

By Lemma \ref {le:3.5}, we have $m(G)=m([T_{G}])+m(C_{q})$. Note that $m=\frac{q}{2}+l$. Since $2m=q+2l$ and the definition of $\mathbf{U_{2m}}$, by the similar proof as in Lemma \ref{le:3.6}, we have $\mathbf{U_{2m}}=F_{1}^{\sigma}\cup F_{3}^{\sigma}$. If $(u,\sigma)\in F_{1}^{\sigma}$, then $p(u)=m,~s(u)=c(u)=0$.
If $(u,\sigma)\in F_{3}^{\sigma}$, then $p(u)=l+1,~c(u)=1$. Let $s(u)=2a+1,~a\in \mathbf{Z^{+}}$ (since $C_{q}^{\sigma}$ has an odd number of negative edges), then

$b_{m}=\sum_{(u,\sigma)\in F_{1}^{\sigma}}(-1)^{m}+\sum_{(u,\sigma)\in F_{3}^{\sigma}}(-1)^{l+1+2a+1}\cdot 2^{1}$ = $(-1)^{m}\cdot|F_{1}^{\sigma}|+(-1)^{l}\cdot2\cdot|F_{3}^{\sigma}|=(-1)^{l}((-1)^{m-l}|F_{1}^{\sigma}|+2\cdot|F_{3}^{\sigma}|)$. Since $m-l=\frac{q}{2}$ and $q\equiv 0(mod~4)$, we have $m-l\equiv 0(mod~2)$. Thus $b_{m}=(-1)^{l}(|F_{1}^{\sigma}|+2\cdot|F_{3}^{\sigma}|)\neq 0$.

\textbf{(2).} $m(T_{G})\neq m([T_{G}])$.

Equivalently, $m(T_{G})=m([T_{G}])+1$, then each maximum matching of $T_{G}$ must cover the unique cyclic vertex. By Lemma \ref{le:3.5}, we mainly discuss the following two cases.

\textbf{Case 1.} $m(G)\neq m(C_{q}^{\sigma})+m([T_{G}])$,
that is, $m(G)> m(C_{q}^{\sigma})+m([T_{G}])$,~$m>\frac{q}{2}+l,~2m>q+2l$. Then the basic subgraphs with $2m$ vertices have no cycle. Which shows that $F_{3}^{\sigma}=\phi$ and $\mathbf{U_{2m}}=F_{1}^{\sigma}$. Then $p(u)=m$ and $s(u)=c(u)=0$ for any $(u,\sigma) \in \mathbf{U_{2m}}$. Thus
$b_{m}=\sum_{(u,\sigma)\in F_{1}^{\sigma}}(-1)^{m}=(-1)^{m}|F_{1}^{\sigma}|\neq 0$.

\textbf{Case 2.} $m(G)=m(C_{q}^{\sigma})+m([T_{G}])$ and $M\cap E_{1}\neq \phi,~\exists~M\in F_{1}$ (equivalently, $F_{1}^{'}\neq \phi$).
 Note that $m=l+\frac{q}{2}$, there exists some basic subgraphs with $2m$ vertices, such that contain $C_{q}^{\sigma}$ as a subgraph. Then
$\mathbf{U_{2m}}=F_{1}^{\sigma}\cup F_{3}^{\sigma}$. If $(u,\sigma)\in F_{1}^{\sigma}$, then $p(u)=m, s(u)=c(u)=0$.
If $(u,\sigma)\in F_{3}^{\sigma}$, we have $p(u)=l+1,~c(u)=1$. Let $s(u)=2a+1,~a\in \mathbf{Z^{+}}$ (since $C_{q}^{\sigma}$ has an odd number of negative edges). Then

$b_{m}=\sum_{(u,\sigma)\in F_{1}^{\sigma}}(-1)^{m}+\sum_{(u,\sigma)\in F_{3}^{\sigma}}(-1)^{l+1+2a+1}\cdot 2^{1}$ = $(-1)^{m}\cdot|F_{1}^{\sigma}|+(-1)^{l}\cdot2\cdot|F_{3}^{\sigma}|=(-1)^{l}((-1)^{m-l}|F_{1}^{\sigma}|+2\cdot|F_{3}^{\sigma}|)$. Since $q$ is even, we mainly discuss the following two subcases.

\textbf{Subcase 2.1.} If $q\equiv 0(mod~4)$, then $m-l=\frac{q}{2}\equiv 0(mod~2)$. Thus $$b_{m}=(-1)^{l}\cdot(|F_{1}^{\sigma}|+2\cdot|F_{3}^{\sigma}|)\neq 0$$

\textbf{Subcase 2.2.} If $q\equiv 2(mod~4)$, then $m-l=\frac{q}{2}\equiv 1(mod~2)$. Thus

$b_{m}=(-1)^{l+1}(|F_{1}^{\sigma}|-2\cdot|F_{3}^{\sigma}|)=(-1)^{l+1}(|F_{1}^{\sigma}|-2\cdot|F_{2}^{\sigma}|)$, $|F_{3}|=|F_{2}|$. Note that $F_{1}^{'} \neq \phi$, by Corollary \ref{co:1}, we have $|F_{1}|>2|F_{2}|$. Hence $b_{m}\neq 0$,

Based on the above conclusions, we have $b_{m}\neq 0$. Thus $\eta(G,\sigma)=|V(G)|-2m(G)$.\quad $\square$

By Lemmas \ref{le:3.2}, \ref{le:3.3}, \ref{le:3.6} and \ref{le:3.7}, we obtained Theorem \ref{thm:1.1}.

\section{Proof of Theorem 1.2}

In this section, we mainly investigate the unbalanced signed graphs. In order to proof Theorem \ref{thm:1.2}, we need to introduce some definitions and lemmas.

Let $G$ be the graph with some pendant vertices, and $G$ contains at least a cycle. For any pendant vertex $u$, we will give definitions of two types of pendant vertex $u$:
\begin{defn}\label{de:2}

\begin{enumerate}[(1)]
  \item If the quasi-pendant vertex of $u$ does not lie on the cycle, then $u$ is of Type \uppercase\expandafter{\romannumeral1}.
  \item If the quasi-pendant vertex of $u$ lies on the cycle, then $u$ is of Type \uppercase\expandafter{\romannumeral2}.
\end{enumerate}
\end{defn}

\begin{lem} \label{le:4.1}
Let $(G,\sigma)$ be an unbalanced signed graph with a pendant vertex $u$, and $v$ is adjacent to $u$. If $u$ is of Type \uppercase\expandafter{\romannumeral1}, then
$\eta(G,\sigma)=|V(G)|-2m(G)+2c(G)-s$ if and only if $\eta (G-u-v,\sigma)=|V(G-u-v)|-2m(G-u-v)+2c(G-u-v)-s$. $(0\leq s\leq 3c(G))$.
\end{lem}

\textbf{Proof.} By Lemmas \ref {le:2.8} and \ref {le:2.12},
$$|V(G)|=|V(G-u-v)|+2,$$
$$m(G)=m(G-u-v)+1,$$
$$c(G)=c(G-u-v).$$

By Lemma \ref {le:2.4},

$\eta (G-u-v,\sigma)=\eta (G,\sigma)=|V(G)|-2m(G)+2c(G)-s=|V(G-u-v)|+2-2(m(G-u-v)+1)+2c(G-u-v)-s=|V(G-u-v)|-2m(G-u-v)+2c(G-u-v)-s$.\quad $\square$

\begin{lem} \label{le:4.2}
Let $(G,\sigma)$ be an unbalanced signed graph with a pendant vertex $u$, and $v$ is adjacent to $u$. If $u$ is of Type \uppercase\expandafter{\romannumeral2}, then
$\eta(G,\sigma)\leq |V(G)|-2m(G)+2c(G)-2$.
\end{lem}

\textbf{Proof.} We prove the result by contradiction. From Lemma \ref{le:2.2}, suppose on the contrary, there exists some unbalanced signed graphs $(H,\sigma)$ with the nullity $\eta(H,\sigma) =|V(H)|-2m(H)+2c(H)-s$, $s=0,1$. Since $v$ lies on the cycle of $(H,\sigma)$, by Lemmas \ref {le:2.8} and \ref {le:2.12}, we have
$$|V(H)|=|V(H-u-v)|+2,$$
$$m(H)=m(H-u-v)+1,$$
$$c(H)\geq c(H-u-v)+1.$$
By Lemma \ref {le:2.4}, we have $\eta (H-u-v,\sigma)=\eta (H,\sigma)=|V(H)|-2m(H)+2c(H)-s\geq |V(H-u-v)|+2-2(m(H-u-v)+1)+2c(H-u-v)+2-s\geq |V(H-u-v)|-2m(H-u-v)+2c(H-u-v)+1$.
Which contradicts Lemma \ref{le:2.2}. Therefore, for any unbalanced signed graph $(G,\sigma)$ with a pendant vertex $u$  of Type \uppercase\expandafter{\romannumeral2}, then
$\eta(G,\sigma)\leq |V(G)|-2m(G)+2c(G)-2$.\quad $\square$

\begin{lem} \label{le:4.3}
Let $(G,\sigma)$ be an unbalanced signed graph without pendant vertices. If $\eta(G,\sigma)\neq |V(G)|-2m(G)+2c(G)$,  $c(G)\geq2$, then there exists a vertex $x$ on the cycle in $(G,\sigma)$ such that $\eta(G-x,\sigma)\neq |V(G-x)|-2m(G-x)+2c(G-x)$.
\end{lem}

\textbf{Proof.}
\textbf{(1).} $g(G)=3$,~$c(G)\geq 2$.

Denote by $C_{q}^{\sigma}$ the cycle of $(G,\sigma)$ and $q=3$. Since $c(G)\geq2$, there exists a vertex $x$ on another cycle in $(G,\sigma)$ such that $C_{q}^{\sigma}$ is a subgraph of $(G-x,\sigma)$. Which shows that $(G-x,\sigma)$ does not satisfy in Lemma \ref{le:2.3}(2), then $$\eta(G-x,\sigma)\neq |V(G-x)|-2m(G-x)+2c(G-x).$$

\textbf{(2).} $g(G)\geq4$,~$c(G)\geq 2$.

Since $\eta(G,\sigma)\neq |V(G)|-2m(G)+2c(G)$, $(G,\sigma)$ does not satisfy at least one of the three conditions in Lemma \ref{le:2.3}. 

\textbf{Case 1.} If $(G,\sigma)$ does not satisfy Lemma \ref{le:2.3}(1).

 Note that $(G,\sigma)$ contains at least two vertex-joint cycles, called $C_{k}^{\sigma}$ and $C_{s}^{\sigma}$ $(k,s\geq4)$. Let $G[C_{k}^{\sigma},C_{s}^{\sigma}]$ be the graph induced by $V(C_{k}^{\sigma})$ and $V(C_{s}^{\sigma})$.

 \textbf{Subcase 1.1.} $c(G)=2$.

 Since $(G,\sigma)$ contains no pendant vertices, then $(G,\sigma)$ is the union of an $\infty$-graph (or a $\theta$-graph) and some isolated vertices (if any), which implies that $G[C_{k}^{\sigma},C_{s}^{\sigma}]$ is either $\infty^{\sigma}(p,1,q)$ or $\theta^{\sigma}(p,l,q)$.

\textbf{Subcase 1.1.1.} $G[C_{k}^{\sigma},C_{s}^{\sigma}]=\theta^{\sigma}(1,1,1)$.

Since $(G,\sigma)$ is unbalanced, $(G,\sigma)$ contains at least an unbalanced cycle, says $C_{q}^{\sigma}$, and $q=4$. Let $x$ be a vertex of another cycle in $(G,\sigma)$, $x\notin V(C_{q}^{\sigma})$. Then $C_{q}^{\sigma}$ is a subgraph of $(G-x,\sigma)$. Which shows that $(G-x,\sigma)$ does not satisfy Lemma \ref{le:2.3}(2), then $$\eta(G-x,\sigma)\neq |V(G-x)|-2m(G-x)+2c(G-x).$$

\textbf{Subcase 1.1.2.} $G[C_{k}^{\sigma},C_{s}^{\sigma}] \neq\theta^{\sigma}(1,1,1)$.

As shown in Figure 1, there exists a vertex $x$ on the cycle in $(G,\sigma)$, such that $(G-x,\sigma)$ contains a pendant vertex with Type \uppercase\expandafter{\romannumeral2}. By Lemma \ref{le:4.2}, we can get that $$\eta(G-x,\sigma)\neq |V(G-x)|-2m(G-x)+2c(G-x).$$

\textbf{Subcase 1.2.} $c(G)\geq 3$.

\textbf{Subcase 1.2.1.} In $(G,\sigma)$, there exists a vertex $x$ on the cycle such that $x\notin G[C_{k}^{\sigma},C_{s}^{\sigma}]$.

 Which implies that $G[C_{k}^{\sigma},C_{s}^{\sigma}]$ is a subgraph of $(G-x,\sigma)$. Then $(G-x,\sigma)$ does not satisfy Lemma \ref{le:2.3}(1), hence $$\eta(G-x,\sigma)\neq |V(G-x)|-2m(G-x)+2c(G-x).$$
For example, as shown in Figure 2, the signed graphs with $T_{i}(i=1,\ldots,4)$ as a underlying graph, which contains a vertex $x$ on the cycle and $x\notin G[C_{k}^{\sigma},C_{s}^{\sigma}]$.

\textbf{Subcase 1.2.2.} In $(G,\sigma)$, let $x$ be a vertex on any cycle, then $x\in G[C_{k}^{\sigma},C_{s}^{\sigma}]$.

Which implies that any cycle of $(G,\sigma)$ is the subgraph of $G[C_{k}^{\sigma},C_{s}^{\sigma}]$. Since $c(G)\geq3$, $(G,\sigma)$ contains one of eight types of bases of tricyclic graphs as a underlying subgraph. As shown in Figure 2, the graph $T_{j}$ can be regard as two vertex-joint cycles, where $j=5,\ldots,8$. Which implies that the tricyclic graphs $T_{j}$ is a underlying subgraph of $G[C_{k}^{\sigma},C_{s}^{\sigma}]$. As shown in Figure 2, there exists a 
vertex $x$ of $T_{j}$ such that $T_{j}-x$ also contains two vertex-joint cycles. Hence there exists a vertex $x$ on the cycle of $(G,\sigma)$ such that $(G-x,\sigma)$ does not satisfy Lemma \ref{le:2.3}(1), then $$\eta(G-x,\sigma)\neq |V(G-x)|-2m(G-x)+2c(G-x).$$

\textbf{Case 2.} $(G,\sigma)$ satisfies Lemma \ref{le:2.3}(1) but does not satisfy Lemma \ref{le:2.3}(2).

\textbf{Subcase 2.1.} There exists at least a balanced cycle in $(G,\sigma)$.

\textbf{Subcase 2.1.1.} There exists a balanced cycle in $(G,\sigma)$, says $C_{p}^{\sigma}$ and $p \neq0(mod~4)$. Let $x$ be a vertex on another cycle and $x\notin V(C_{p}^{\sigma})$ (since all cycles are pairwise vertex-disjoint), then $C_{p}^{\sigma}$ is the subgraph of $(G-x,\sigma)$. Which shows that $(G-x,\sigma)$ does not satisfy Lemma \ref{le:2.3}(2), then$$\eta(G-x,\sigma)\neq |V(G-x)|-2m(G-x)+2c(G-x).$$

\textbf{Subcase 2.1.2.} The length of all balanced graph is the multiple of $4$. Then $(G,\sigma)$ contains at least an unbalanced cycle in $(G,\sigma)$, says $C_{q}^{\sigma}$ and $q \neq2(mod~4)$. Let $x$ be a vertex of the balanced cycle and $x\notin V(C_{q}^{\sigma})$ (Since all cycles are pairwise vertex-disjoint), then $C_{q}^{\sigma}$ is the subgraph of $(G-x,\sigma)$. Which shows that $(G-x,\sigma)$ does not satisfy Lemma \ref{le:2.3}(2), then $$\eta(G-x,\sigma)\neq |V(G-x)|-2m(G-x)+2c(G-x).$$

\textbf{Case 2.2.} Cycles of $(G,\sigma)$ are unbalanced.

$(G,\sigma)$ contains at least an unbalanced cycle, says $C_{q}^{\sigma}$ and $q \neq2(mod~4)$. Let $x$ be a vertex of another unbalanced cycle in $(G,\sigma)$ and $x\notin V(C_{q}^{\sigma})$ (all cycles are pairwise vertex-disjoint). Since $c(G)\geq2$, then $C_{q}^{\sigma}$ is the subgraph of $(G-x,\sigma)$. Which shows that $(G-x,\sigma)$ does not satisfy Lemma \ref{le:2.3}(2), then $$\eta(G-x,\sigma)\neq |V(G-x)|-2m(G-x)+2c(G-x).$$

\textbf{Case 3.} $(G,\sigma)$ satisfies (1) and (2) of Lemma \ref{le:2.3}, but dose not satisfy (3) of Lemma \ref{le:2.3}.

Note that $m(T_{G})\neq m([T_{G}])$, that is, $m(T_{G})\geq m([T_{G}])+1$.

If $E(T_{G})=\phi$, then $(G,\sigma)$ is the union of some vertex-disjoint signed cycles and isolated vertices. Hence $m(T_{G})=m([T_{G}])=0$, a contradiction. Therefore, we only consider $E(T_{G})\neq\phi$. In $T_{G}$, each maximum matching must cover at least a pendant vertex. Otherwise, there exists an $M$-augmenting path in $G$, which contradicts Lemma \ref{le:2.13}. Let $u$ be a pendant vertex of $T_{G}$. Since $(G,\sigma)$ contains no pendant vertices, we have $u \in W_{G}$. Suppose that $C_{q}^{\sigma}$ the pendant cycle of $(G,\sigma)$ is correspond to the vertex $u$ of $T_{G}$. Let $v$ be the unique vertex with degree three in $C_{q}^{\sigma}$, and let $x$ be a vertex on the cycle $C_{q}^{\sigma}$. Then $T_{G-x}$ is the graph obtained from $T_{G}$ and $C_{q}^{\sigma}-x$ by identifying $u$ and $v$ as a vertex, we mainly discuss the following two subcases.

\textbf{Subcase 3.1.} Each maximum matching of $T_{G}$  cover all pendant vertices.

Let $x$ be a vertex on the cycle $C_{q}^{\sigma}$, and $x$ is adjacent to $v$. Noted that $C_{q}^{\sigma}$ is an even cycle, then $C_{q}^{\sigma}-v-x$ is an odd path and has a perfect matching. By the definition of $T_{G}$, which shows that the maximum matching of $T_{G-x}$ is the union of the maximum matching of $T_{G}$ and the maximum matching of $C_{q}^{\sigma}-v-x$. Then $$m(T_{G-x})=m(T_{G})+m(C_{q}^{\sigma}-v-x).$$ Hence each maximum matching of $T_{G-x}$ must cover some vertices in $W_{G-x}$, 
we have $m(T_{G-x})\neq m([T_{G-x}])$. Which shows that $(G-x,\sigma)$ does not satisfy Lemma \ref{le:2.3}(3), then $$\eta(G-x,\sigma)\neq |V(G-x)|-2m(G-x)+2c(G-x).$$

\textbf{Subcase 3.2.} There exists some maximum matchings of $T_{G}$, denoted by $M_{i}(T_{G})$ $(i=1,\ldots,r)$, such that the pendant edge $wv\notin M_{i}(T_{G})$, where $r\geq 1$ and $v$ is a pendant vertex of $T_{G}$.

Let $y$ be a vertex on the cycle $C_{p}^{\sigma}$ and $d(v_{0},y)=2$, where $C_{p}^{\sigma}$ the cycle of $(G,\sigma)$ is correspond to $v$ the vertex of $T_{G}$, and $v_{0}$ is the unique vertex with degree three in $C_{p}^{\sigma}$. By the definition of $T_{G-y}$, which shows that the maximum matching of $T_{G-y}$ is the union of $M_{i}(T_{G})$ $(i\in \{1,\ldots,r\})$ and the maximum matching of $C_{p}^{\sigma}-y$. Then $$m(T_{G-y})=m(T_{G})+m(C_{p}^{\sigma}-y).$$ For any $i\in \{1,\ldots,r\}$, then $M_{i}(T_{G})$ must cover some vertices in $W_{G-y}$ by Lemma \ref{le:2.13}. Which implies that each maximum matching of $T_{G-y}$ must cover some vertices in $W_{G-y}$, we have $m(T_{G-y})\neq m([T_{G-y}])$. Which shows that $(G-y,\sigma)$ does not satisfy Lemma \ref{le:2.3}(3), then $\eta(G-y,\sigma)\neq |V(G-y)|-2m(G-y)+2c(G-y)$. \quad $\square$

\begin{lem} \label{le:4.4}
Let $(G,\sigma)$ be an unbalanced signed graph without pendant vertices. Then $$\eta(G,\sigma)\neq |V(G)|-2m(G)+2c(G)-1.$$
\end{lem}

\textbf{Proof.} We apply induction on $c(G)$ to prove the result.

If $c(G)=1$, then $G=C_{k}\cup (n-k)K_{1}$~$(3\leq k\leq n)$. By Theorem \ref{thm:1.1}, we have $$\eta(G,\sigma)\neq |V(G)|-2m(G)+2c(G)-1.$$

If $c(G)\geq 2$, assume that the conclusion is true for $c(G)\leq k$. Next, we need to prove the conclusion is true for $c(G)= k+1$. Suppose on the contrary, there exists some unbalanced signed graphs $(H,\sigma)$ with $c(H)=k+1$ such that $\eta(H,\sigma)=|V(H)|-2m(H)+2c(H)-1$.

Let $x$ be any vertex on the cycle of $(H,\sigma)$. For the signed graph $(H-x,\sigma)$, by Lemmas \ref{le:2.9} and \ref{le:2.12}, we have
$$|V(H)|=|V(H-x)|+1,$$
$$m(H)\leq m(H-x)+1,$$
$$c(H)\geq c(H-x)+1.$$
By Lemma \ref{le:2.1}, we have $\eta(H-x,\sigma)+1\geq \eta(H,\sigma)=|V(H)|-2m(H)+2c(H)-1 \geq|V(H-x)|+1-2(m(H-x)+1)+2(c(H-x)+1)-1=|V(H-x)|-2m(H-x)+2c(H-x)$.
By Lemma \ref{le:2.2}, then $$\eta(H-x,\sigma)=|V(H-x)|-2m(H-x)+2c(H-x)-s,~s=0,1$$.

If $(H-x,\sigma)$ is balanced, by Lemmas \ref{le:2.14} and \ref {le:2.5}, we have
$\eta(H-x,\sigma) \neq|V(H-x)|-2m(H-x)+2c(H-x)-1$.
We mainly consider that $(H-x,\sigma)$ is unbalanced in the following cases.

\textbf{Case 1.} $(H-x,\sigma)$ contains no pendant vertices.

In this case, $(H-x,\sigma)$ is also unbalanced and $c(H-x)\leq c(H)-1=k$. Applying the induction hypothesis to $(H-x,\sigma)$,
then $\eta(H-x,\sigma)\neq|V(H-x)|-2m(H-x)+2c(H-x)-1$.

\textbf{Case 2.} $(H-x,\sigma)$ contains some pendant vertices.

\textbf{Subcase 2.1.} $(H-x,\sigma)$ contains at least a pendant vertex of Type \uppercase\expandafter{\romannumeral2}.

By Lemma \ref{le:4.2}, we have $\eta(H-x,\sigma) \leq|V(H-x)|-2m(H-x)+2c(H-x)-2$.

\textbf{Subcase 2.2.} All the pendant vertices of $(H-x,\sigma)$ are of Type \uppercase\expandafter{\romannumeral1}.

Suppose that $(H-x,\sigma)$ contains $p$ pendant vertices. For the pendant vertices of Type \uppercase\expandafter{\romannumeral1}, by using Lemma \ref{le:2.4} repeatedly, after $p$ steps,
we obtain a subgraph $(H_{1},\sigma)$ of $(H,\sigma)$. If $(H_{1},\sigma)$ contains no pendant vertices or at least a pendant vertex of Type \uppercase\expandafter{\romannumeral2}, then $(H_{1},\sigma)$ is as we required and we are done. Otherwise, in $(H_{1},\sigma)$, for the pendant vertex of Type \uppercase\expandafter{\romannumeral1}, we continue to use Lemma \ref{le:2.4} repeatedly, we obtain a subgraph $(H_{2},\sigma)$ of $(H_{1},\sigma)$. If $(H_{2},\sigma)$ contains no pendant vertices or at least a pendant vertex of Type \uppercase\expandafter{\romannumeral2}, then $(H_{2},\sigma)$ is as we required and we are done. Otherwise, repeating the above steps until we obtain an unbalanced sign graph $(H_{0},\sigma)$ that meets the requirements.

\textbf{(a).} $(H_{0},\sigma)$ contains at least a pendant vertex of Type \uppercase\expandafter{\romannumeral2}.

By Lemma \ref{le:4.2}, we have $$\eta(H_{0},\sigma)\leq |V(H_{0})|-2m(H_{0})+2c(H_{0})-2.$$  Next, by Lemma \ref{le:4.1},

$$\eta(H_{0},\sigma)\leq |V(H_{0})|-2m(H_{0})+2c(H_{0})-2$$
if and only if $$\eta(H-x,\sigma)\leq |V(H-x)|-2m(H-x)+2c(H-x)-2.$$

\textbf{(b).} $(H_{0},\sigma)$ contains no pendant vertices.

Since $c(H_{0})=c(H-x)\leq c(H)-1=k$ and $H_{0}^{\sigma}$ is unbalanced, applying the induction hypothesis to $H_{0}^{\sigma}$, then $\eta(H_{0},\sigma)\neq |V(H_{0})|-2m(H_{0})+2c(H_{0})-1$.

By Lemma \ref{le:4.1},
$$\eta(H_{0},\sigma)\neq |V(H_{0})|-2m(H_{0})+2c(H_{0})-1$$ if and only if $$\eta(H-x,\sigma)\neq |V(H-x)|-2m(H-x)+2c(H-x)-1.$$

Based on the above results. Let $x$ be any vertex on the cycle of $(H,\sigma)$. For the nullity $\eta(H-x,\sigma)$ of signed graph $(H-x,\sigma)$, either $\eta(H-x,\sigma)\neq |V(H-x)|-2m(H-x)+2c(H-x)-1$ or  $\eta(H-x,\sigma)\leq |V(H-x)|-2m(H-x)+2c(H-x)-2$.

If $\eta(H-x,\sigma)\neq |V(H-x)|-2m(H-x)+2c(H-x)-1$, then $\eta(H-x,\sigma)= |V(H-x)|-2m(H-x)+2c(H-x)$.
On the other hand, since $\eta(H,\sigma)\neq |V(H)|-2m(H)+2c(H)$, by Lemma \ref{le:4.3}, there exists a vertex $y$ on the cycle of $(H,\sigma)$ such that $\eta(H-x,\sigma)\neq |V(H-x)|-2m(H-x)+2c(H-x)$, a contradiction.

If $\eta(H-x,\sigma)\leq |V(H-x)|-2m(H-x)+2c(H-x)-2$, which contradicts equation $\eta(H-x,\sigma)=|V(H-x)|-2m(H-x)+2c(H-x)-s,~s=0,1$.

Therefore, for any unbalanced signed graph $(G,\sigma)$ without pendant vertices, $\eta(G,\sigma)\neq |V(G)|-2m(G)+2c(G)-1$.\quad $\square$


\textbf{The Proof of Theorem \ref{thm:1.2}:} Noted that the unbalanced signed graph $(G, \sigma)$ contains at least an unbalanced cycle, and $c(G)\geq1$. 

\textbf{Case 1.} $(G,\sigma)$ has no pendant vertices, we can obtained the result by Lemma \ref{le:4.4}.

\textbf{Case 2.} $(G,\sigma)$ has some pendant vertices.

By the definitions of two types of pendant vertices, we mainly discuss the following two subcases.

\textbf{Subcase 2.1.} All pendant vertices are of Type \uppercase\expandafter{\romannumeral1}.

Suppose that $(G,\sigma)$ contains $p$ pendant vertices. For the pendant vertex of Type \uppercase\expandafter{\romannumeral1}, by using Lemma \ref{le:2.4} repeatedly, after $p$ steps, we obtain a subgraph $(G_{1},\sigma)$ of $(G,\sigma)$. If $(G_{1},\sigma)$ contains no pendant vertices or at least a pendant vertex of Type \uppercase\expandafter{\romannumeral2}, then $(G_{1},\sigma)$ is as we required and we are done. Otherwise, in $(G_{1},\sigma)$, for the pendant vertex of Type \uppercase\expandafter{\romannumeral1}, we continue to use Lemma \ref{le:2.4} repeatedly, we obtain a subgraph $(G_{2},\sigma)$ of $(G_{1},\sigma)$. If $(G_{2},\sigma)$ contains no pendant vertices or at least a pendant vertex of Type \uppercase\expandafter{\romannumeral2}, then $(G_{2},\sigma)$ is as we required and we are done. Otherwise, repeating the above steps until we obtain an unbalanced sign graph $(G_{0},\sigma)$ that meets the requirements.

\textbf{(a).} $(G_{0},\sigma)$ contains no pendant vertices.

By Lemma \ref{le:4.4}, we have $\eta(G_{0},\sigma)\neq |V(G_{0})|-2m(G_{0})+2c(G_{0})-1$.

\textbf{(b).} $(G_{0},\sigma)$ contains at least a pendant vertex of Type \uppercase\expandafter{\romannumeral2}.

By Lemma \ref{le:4.2}, we have $$\eta(G_{0},\sigma)\leq |V(G_{0})|-2m(G_{0})+2c(G_{0})-2.$$ Hence $$\eta(G_{0},\sigma)\neq |V(G_{0})|-2m(G_{0})+2c(G_{0})-1.$$ Next, by Lemma \ref{le:4.1}, we have

$\eta(G_{0},\sigma)\neq |V(G_{0})|-2m(G_{0})+2c(G_{0})-1$ if and only if $\eta(G,\sigma)\neq |V(G)|-2m(G)+2c(G)-1$.

\textbf{Subcase 2.2.} There exists at least a pendant vertex of Type \uppercase\expandafter{\romannumeral2}.

By Lemma \ref{le:4.2}, we have $\eta(G,\sigma)\neq|V(G)|-2m(G)+2c(G)-1$.\quad $\square$

Combining with Lemmas \ref{le:2.14}, \ref{le:2.5} and Theorem \ref{thm:1.2}, we can obtained the following corollary.

\begin{cor}\label{co:3}
For any signed graph $(G,\sigma)$, then $$\eta(G,\sigma)\neq|V(G)|-2m(G)+2c(G)-1.$$
\end{cor}


\section{Proof of Theorem 1.3}

In order to  proof Theorem \ref{thm:1.3}, we need to introduce some definitions.

Let $K_{1,c+1}^{\sigma}$ be a signed star, $x$ be the center vertex of $K_{1,c+1}^{\sigma}$ and $y_{0},\ldots, y_{c}$ be pendant vertices of $K_{1,c+1}^{\sigma}$,~respectively. 
$O_{1}^{\sigma},\ldots,O_{l_{1}}^{\sigma}$ be a balanced cycle of size 3, $O_{l_{1}+1}^{\sigma},\ldots,O_{l_{1}^{'}}^{\sigma}$ be an unbalanced cycle of size 6. Let $O_{l_{1}^{'}+1}^{\sigma},\ldots,O_{c}^{\sigma}$ be the signed graph with the pendant vertex of Type \uppercase\expandafter{\romannumeral2}, which obtained from $C_{4}^{+}$ a balanced cycle of size 4 by attaching $K_{2}^{\sigma}$.
~Next, we construct a new signed graph $G^{\sigma}$ by identifying $y_{i}$ with a vertex of $O_{i}^{\sigma}$, $y_{j}$ with a pendant vertex of $O_{j}^{\sigma}$, where $i=1,\ldots,l_{1}^{'},~j=l_{1}^{'}+1,\ldots,c$.


\textbf{The Proof of Theorem \ref{thm:1.3}:} According to the definition of $G^{\sigma}$, let $l_{2}=l_{1}^{'}-l_{1},~l_{3}=c-l_{1}^{'}$. In this case, there exists $l_{1}$ balanced cycles $C_{3}^{+}$, $l_{2}$ unbalanced cycles $C_{6}^{-}$ and $l_{3}$ balanced graph $F^{\sigma}$, which is obtained from $C_{4}^{+}$ a balanced cycles of size 4 by attaching $K_{2}^{\sigma}$. Which implies that $y_{0}$ is the unique pendant vertex of $G^{\sigma}$.
Then
$$|V(G^{\sigma})|=3l_{1}+6l_{2}+5l_{3}+2,$$
$$m(G^{\sigma})=l_{1}+3l_{2}+2l_{3}+1,$$
$$c(G^{\sigma})=l_{1}+l_{2}+l_{3}.$$

By Lemmas \ref{le:2.4} and \ref{le:2.11}, we obtain that $\eta(G^{\sigma})=\eta(G^{\sigma}-x-y_{0}))=2l_{2}+l_{3}=|V(G^{\sigma})|
-2m(G^{\sigma})+2c(G^{\sigma})-(3l_{1}+2l_{3})$. Since $l_{i}\geq 0$, $3l_{1}+2l_{3}$ can take over any integer between 0 and $3c(G)$ expect for 1.  \quad $\square$

Combining with Lemma \ref{le:2.5} and Theorem \ref{thm:1.3}, we can obtained the following corollary.
\begin{thm}\label{th:5.1}
For any signed graph $(G,\sigma)$, if $c(G)$ is fixed, then there exists infinitely connected signed graphs, such that $\eta(G,\sigma)=|V(G)|-2m(G)+2c(G)-s$, where $0\leq s\leq 3c(G)$, $s\neq1$.

\end{thm}



\end{document}